\documentclass{article}
\usepackage[utf8]{inputenc}
\usepackage{todonotes}
\usepackage{booktabs}
\usepackage{hyperref}
\usepackage[font=small]{caption}
\usepackage[margin=1in]{geometry}
\usepackage{amsmath}
\usepackage{amsthm}
\usepackage{amssymb}
\usepackage{bbm}
\usepackage{url}
\usepackage{dsfont}
\usepackage{graphicx}
\usepackage{authblk}
\usepackage{subcaption}
\hypersetup{bookmarksopen,bookmarksnumbered,
pdfpagemode=UseOutlines,
colorlinks=true,
urlcolor=blue,
linkcolor=blue,
anchorcolor=blue,
citecolor=blue,
filecolor=blue,
menucolor=blue
}
\usepackage{accents}

\usepackage{cleveref}

\title{Estimating Driver Response Rates to Variable Message Signage at Seattle-Tacoma International Airport}

\author[1]{Soumya Vasisht}
\author[2]{Shushman Choudhury}
\author[1]{Nawaf Nazir}
\author[2]{Stephen Zoepf}
\author[1]{Chase P. Dowling}
\affil[1]{Pacific Northwest National Laboratory\\
          Richland, WA \\
          \{\tt{soumya.vasisht}, \tt{nawaf.nazir}, \tt{chase.dowling}\}\tt{@pnnl.gov}}
\affil[2]{Lacuna Technologies, Inc. \{\tt{shushman.choudhury}, \tt{stephen.zoepf}\}\tt{@lacuna.ai}}

\date{}

\begin{document}

\maketitle

\begin{abstract}
We apply Bayesian Linear Regression to estimate the response rate of drivers to variable message signs at Seattle-Tacoma International Airport, or SEA. Our approach uses vehicle speed and flow data measured at the entrances of the arrival and departure-ways of the airport terminal, and sign message data. Depending on the time of day, we estimate that between 5.5 and 9.1\% of drivers divert from departures to arrivals when the sign reads ``departures full, use arrivals'', and conversely, between 1.9 and 4.2\% of drivers divert from arrivals to departures. Though we lack counterfactual data (i.e., what would have happened had the diversionary treatment not been active), adopting a causal model that encodes time dependency with prior distributions rate can yield a measurable effect.
\end{abstract}


\section{Questions}

Major airports devote significant resources to mitigate congestion~\cite{harris2017mesoscopic}. Their interventions lie on a spectrum of time-scales: bespoke responses to individual incidents (e.g., dispatching operations personnel) at one end and long-term strategic initiatives (e.g., building new infrastructure) at the other~\cite{ugirumurera2021modeling}. In between the two are a range of tactical decisions that may rely on both proactive and reactive actions (e.g., constructing a permanent digital signboard that displays messages in response to real-time changes). We focus on this category of interventions, which currently tend to be based largely on heuristics and domain intuition. Quantifying the effect of such interventions in a numerically robust, reproducible, and interpretable manner is a crucial step towards automated decision-making that can optimize for the objectives that airports care about.

We wish to measure the effectiveness of variable message signage at diverting traffic between arrivals and departures using historical vehicle speed and flow (volume per unit time) and sign data obtained from SEA. Measuring this effectiveness will help implement automated traffic diversion control that optimizes for congestion, sojourn, and/or emissions.
A time delay exists between interventions, i.e., traffic diversion messages, and measurable responses, i.e., vehicle flow entering the roadways. Thus, while we can distinguish between congested and uncongested regimes, naively analyzing the data does not indicate if the system responds. Since we lack counterfactual data, we need to make a prior assumption on 1) the causal relationship between the timing of the sign and the response rate and 2) the distribution of the effectiveness of the signage at diverting drivers.

\begin{figure*}
\centering
    \captionsetup[subfigure]{justification=centering}
    \begin{subfigure}{\textwidth}
        \centering
        \includegraphics[width=0.75\textwidth]{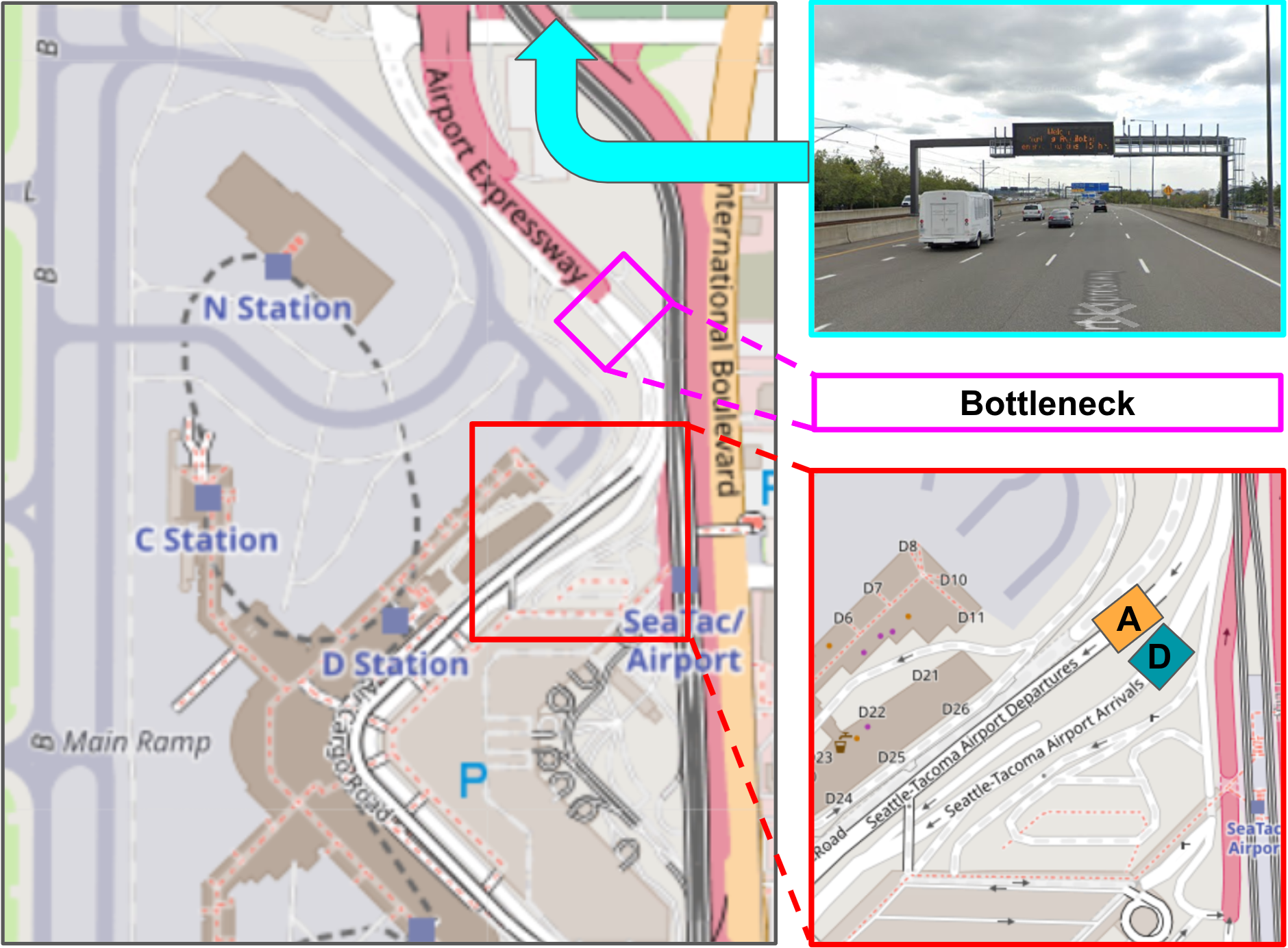}
        \caption{The physical layout of the VMS sign, traffic bottleneck, and arrival and departure drives and the locations of respective speed and flow sensors labeled A and D.}
        \label{fig:seatac-illustr-real}
    \end{subfigure}
    
    \begin{subfigure}{\textwidth}
        \centering
        \includegraphics[width=0.99\textwidth]{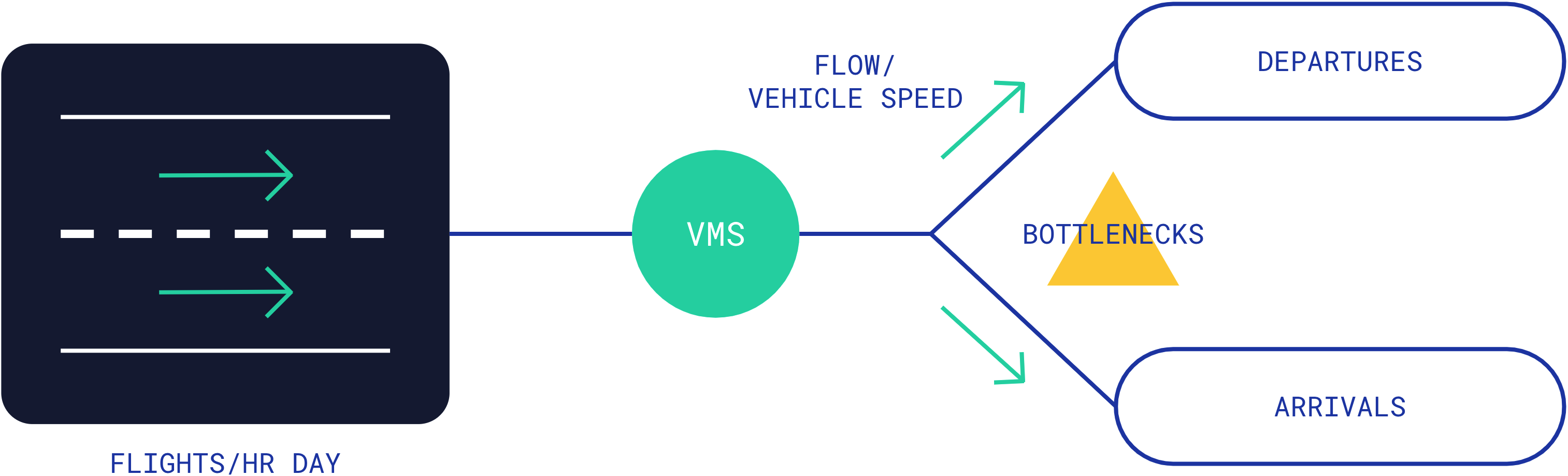}
        \caption{A simplified diagram of the SEA ground traffic flow with the salient components.}
        \label{fig:seatac-illustr-artif}
    \end{subfigure}
    \label{fig:seatac-illustr}
    \caption{}
\end{figure*}

\section{Methods}
The difficulty of our setting stems largely from the available data: vehicle flow and median speed over 4 months, binned to 15-minute intervals, and time-stamped records of messages that request traffic to divert. \textit{We observe no per-vehicle information whatsoever}. Neither do we know the intended and eventual destinations of incoming vehicles, nor can we estimate them. We also cannot run a controlled experiment and must only work with historical observational data.

We frame our problem as estimating the effect of a multi-timestep intervention on a time-varying system. We first quantify the statistical significance of the treatment impact observed on three traffic-related metrics during the intervention period using a two-tailed independent T-test. The first two are the vehicle flow and median speed. We also use a third metric that combines speed and flow as a single measure of congestion. This metric, the so-called \textit{critical ratio}, is the ratio of median speed to the critical speed, which is the speed threshold (given the current flow) below which the system is likely congested~\cite{kerner2009introduction}. The lower the critical ratio than 1, the more congested the system (details omitted for brevity). We posit a null hypothesis about the mean change in these metrics before and during intervention, i.e, the time bins in which the diversionary message remains active. A t-statistic outside the 95\% confidence region (p-value less than 0.05) is evidence of a causal relationship but does not alone yield a function that maps the current system state to an expected number of drivers who will divert if diversion is signaled. 
Therefore we additionally assume a plausible causal model where the treatment is the diversion and the outcome is a function of the difference in incoming flow between Departures and Arrivals. Finally, we use Bayesian regression to estimate the average treatment effect overall and controlling for hour-of-day.

\subsection{Estimating Treatment Effects}

The following is a Bayesian Linear Regression model for estimating time-varying treatment effects:

\begin{equation}
    \begin{aligned}
    \text{Outcome}_t &\sim \mathrm{Normal}\left(\alpha + \beta \cdot \mathrm{Outcome}_{t-1,\ldots} + \gamma \cdot \mathrm{Treatment}_{t, t-1, \ldots}, \sigma \right) \\
    \alpha, \beta, \gamma, \sigma &\sim \mathrm{Priors},
    \end{aligned}
    \label{eq:bayes-reg-general}
\end{equation}
where priors encode domain knowledge~\cite{mcelreath2020statistical}. We have two kinds of treatments, for diverting from Departures to Arrivals (denoted TD) and for the opposite (denoted TA). \textcolor{black}{At most one treatment can be active at a given time}. We also rule out terms beyond time $t-1$. Flow and speed are binned to 15-minute intervals; any driver making a choice in the current time-step would have seen the message in at most the previous time-step.

\begin{table}[t]
    \centering
    \begin{tabular}{@{} lcccrlccc @{}}
        \toprule
        \multicolumn{4}{c}{Dep $\rightarrow$ Arr} && \multicolumn{4}{c}{Arr $\rightarrow$ Dep}\\
        \cmidrule{1-4} \cmidrule{6-9} Variable & t-value & p-value & 95\% CI && Variable & t-value & p-value & 95\% CI\\
        \midrule
        Dep. Speed   &  5.88     & 1.05e-6  & (8.23, 16.9)  && Arr. Speed & 4.59   &  2.15e-5    &  (5.70, 14.50)   \\
        Dep. Flow   &  1.07     &  0.28 & (-10.6, 34.9)   && Arr. Flow & -3.86   &  2.0e-4    &  (-58.2, -18.5)   \\
        Dep. CR   &  5.92    &  8.9e-7 & (0.22, 0.45)  && Arr. CR & 4.30    &  5.96e-05     &  (0.09, 0.26)  \\
        \bottomrule
    \end{tabular}
    \caption{Statistical significance test results on the changes in median speeds, flows and critical ratios (CR) in response to TA and TD. The test numbers indicate that the overall congestion improves significantly in departures and only marginally in arrivals.}
    \label{tab:t-test}
\end{table}

Treatments attempt to redirect traffic between roadways, thus our outcome should depend on the difference in flow, i.e., $\Delta q_t = q^{\mathrm{dep}}_t - q^{\mathrm{arr}}_t$. Note, however, that treatments are deployed when congestion builds up, e.g., TD becomes active when $\Delta q$ is already increasing. Regressing $\Delta q$ on TD \emph{might erroneously conclude that treating Departures sends more traffic towards it}, which is highly implausible.

Instead, we posit that our outcome variable should be the rate of change of $\Delta q$, i.e., $\Delta q^{\prime}_t = \Delta q_t - \Delta q_{t-1}$. Our regression model is then:
\begin{equation}
    \Delta q_t^{\prime} \sim \mathrm{Normal}\left(\alpha + \beta \cdot \Delta q_{t-1}^{\prime} + \gamma_{D} \cdot \mathrm{TD}_{t,t-1} + \gamma_{A} \cdot \mathrm{TA}_{t,t-1}, \sigma \right).
    \label{eq:bayes-reg-spec}
\end{equation}
\textcolor{black}{Here, the posterior distributions of $\gamma_D$ and $\gamma_A$ encode the average estimated effect of the respective treatments, in the current and previous time bin, on the rate of change of $\Delta q_t$. If the treatments work as intended, we expect $\gamma_D < 0$ and $\gamma_A > 0$ (diverting from Departures to Arrivals would reduce the rate-of-change of difference between $q^{\mathrm{dep}}$ and $q^{\mathrm{arr}}$ and vice versa) and $|\gamma_D|, |\gamma_A| > \beta$ (treatments should explain most of the variation when active, not the outcome at the previous time-step). If so, then $-\gamma_D$ and $\gamma_A$ are our estimates of the average number of vehicles that respond to the treatment, i.e., that divert from one facility to the other when the corresponding message is active. Since we do not have counterfactual information from controlled experiments, these regression co-efficients are our best available estimate of the treatment effects.}

Our hypothesis about the first-order rate of change is based on intuition, not a physical scientific model. Therefore, we run another regression with the second-order rate of change $\Delta q^{\prime \prime}_t = \Delta q^{\prime}_t - \Delta q^{\prime}_{t-1}$, to determine if treatments affect second or higher orders. For this regression, we only care about the significance of treatments, i.e. their relative rather than absolute co-efficients. Thus, we standardize the outcome $\Delta q^{\prime \prime}_t \equiv \frac{\Delta q^{\prime \prime}_t - \mathrm{mean}(\Delta q^{\prime \prime}_t)}{\mathrm{stdev}(\Delta q^{\prime \prime}_t)}$, such that $\Delta q^{\prime \prime}_t$ has mean $0$ and standard deviation $1$~\cite{gelman1995bayesian}. All experiments use the PyMC Python library~\cite{patil2010pymc}.

\textcolor{black}{For the outcome variables, we have data for every 15 minute interval for a little over 4 months, i.e., nearly $12000$ time intervals. For the treatments, we have precise time-stamped records of when each message was displayed. Both kinds of treatments vary widely in duration and time of deployment. Most instances lie between $20$ minutes and $3$ hours, but a handful of both kinds last for longer. All Departures treatments are deployed starting between 5 am and 2 pm, and most Arrivals treatments between 8 pm and midnight. Since we know when each treatment starts and ends, there are some 15 minute intervals where treatments are only partially active; we control for this by encoding TA/TD to be the corresponding fraction of 15 minutes, i.e., $\mathrm{TA}/\mathrm{TD} \in [0, 1]$, not $\{0, 1\}$.}

\section{Findings}

\begin{figure}[t]
\centering
\begin{subfigure}[h]{0.95\linewidth}
\includegraphics[width=\linewidth]{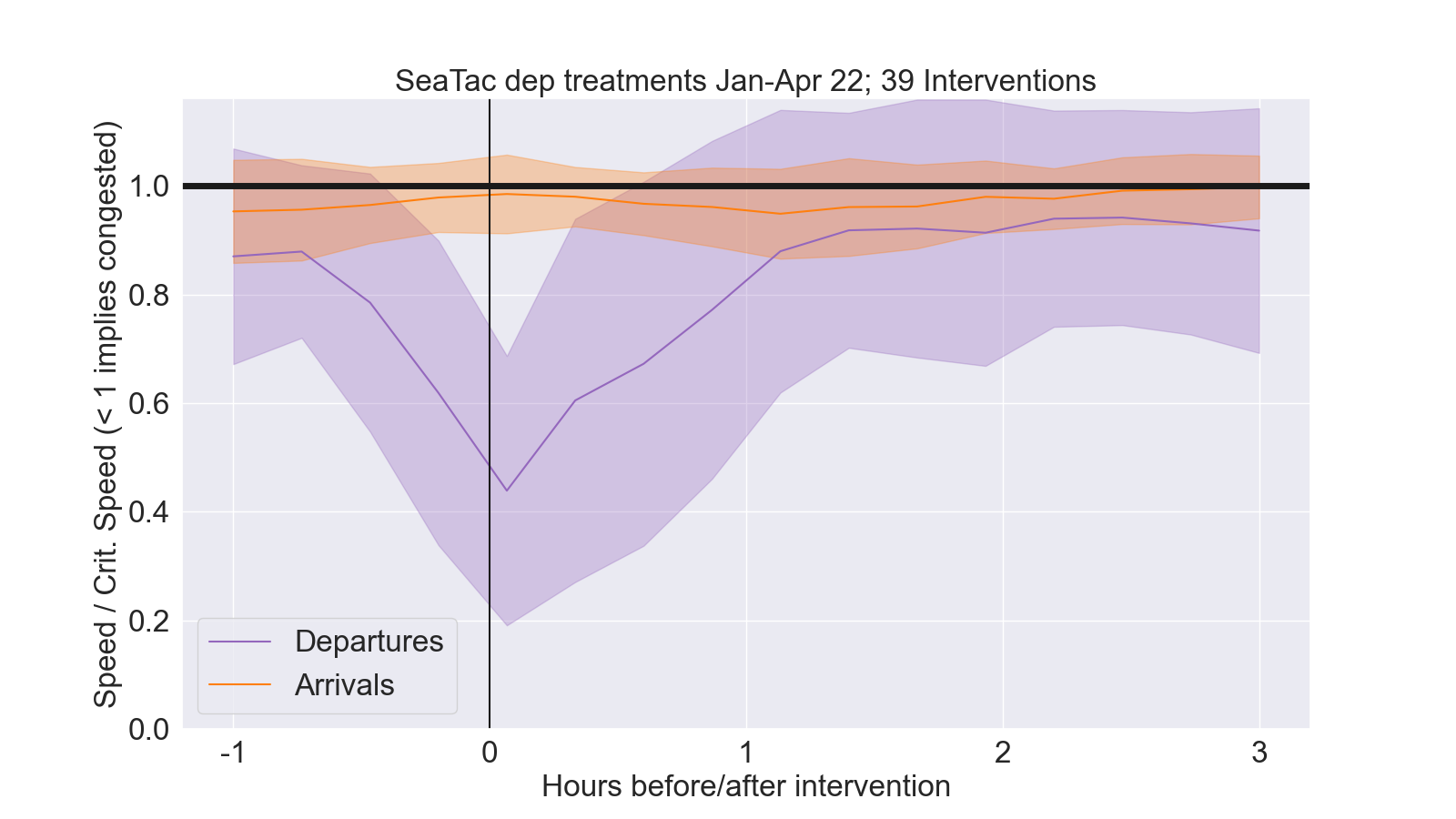}
\caption{Treatment: Departures to Arrivals}
\end{subfigure}
\hfill
\begin{subfigure}[h]{0.95\linewidth}
\includegraphics[width=\linewidth]{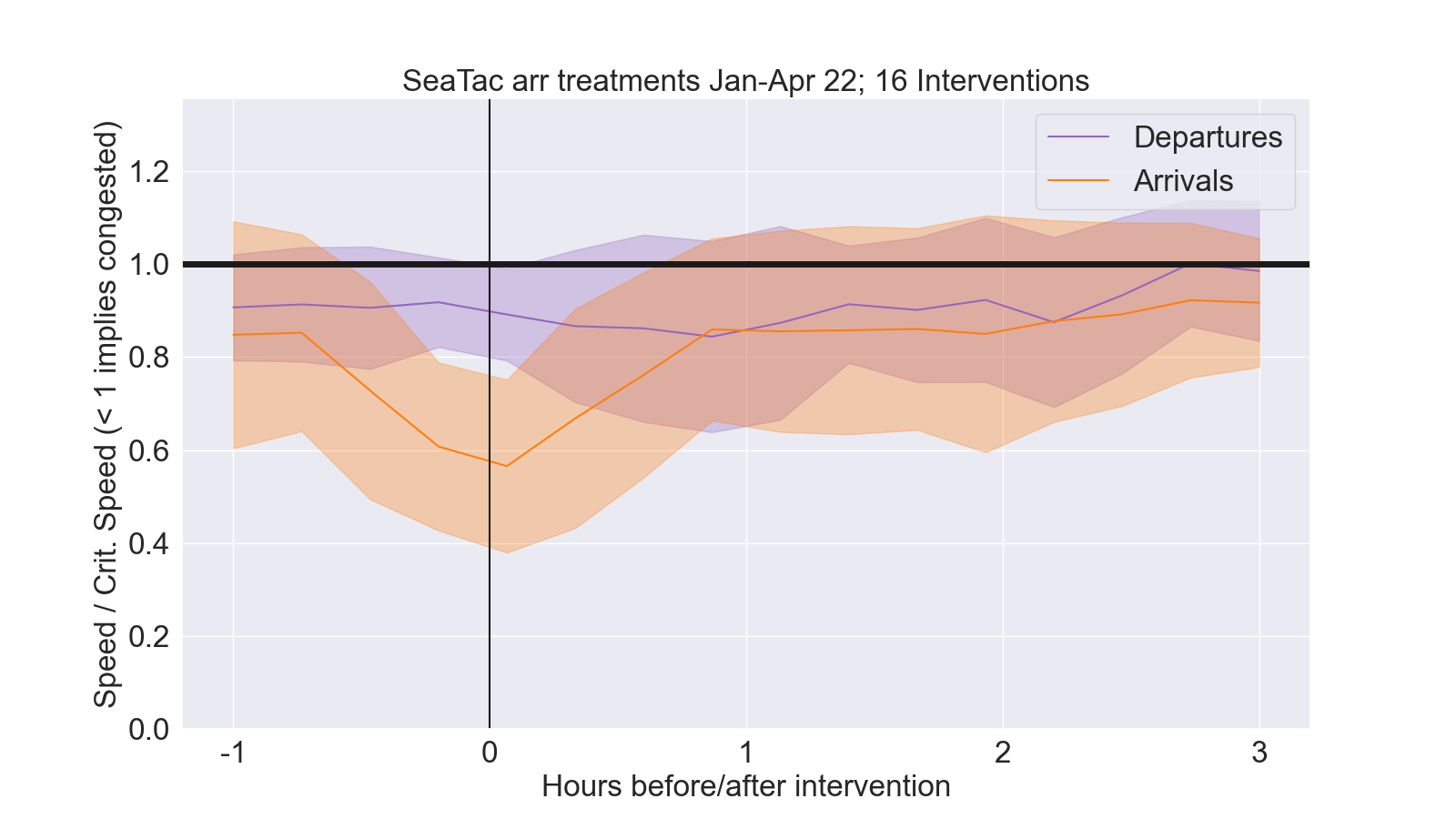}
\caption{Treatment: Arrivals to Departures}
\end{subfigure}%
\caption{The evolution of the speed ratio during the intervention period.}
\label{fig:congestion}
\end{figure}

\subsection{Statistical Significance Test}
\Cref{tab:t-test} summarizes the results of the T-test. The observed effect appears to be consistent for different variables of interest. The changes in median speeds in response to TA and TD are statistically significant, with a low p-value, indicating sufficient evidence to reject the null hypothesis. The departure flow remains unchanged and the arrival flow appears to decrease when the treatment is in effect, which implies some measurable adherence to TD and TA. Furthermore, \textcolor{black}{the $95\%$ confidence intervals (CI) of the difference in means for the critical ratio metric for departures and arrivals demonstrate that overall congestion improves significantly in departures and marginally in arrivals.} This analysis also uncovered some missed opportunities, as noted in \Cref{fig:congestion}. Despite the positive effect of the treatment on the overall system state, the intervention comes in too late, i.e., after the system is already congested. A more proactive, predictive control strategy could decrease the time spent in a congested state, leading to lower overall travel time and many other environmental and logistical benefits.

\begin{table}[t]
    \centering
    \begin{tabular}{@{} lcrcrcclcrcr @{}}
        \toprule
        \multicolumn{5}{c}{Dep $\rightarrow$ Arr} &&& \multicolumn{5}{c}{Arr $\rightarrow$ Dep}\\
        \cmidrule{1-5} \cmidrule{8-12} Hour-of-Day && Inflow/hr && Diverted/hr &&& Hour-of-Day && Inflow/hr && Diverted/hr\\
        \midrule
        0600-0700   &&  925     &&  9.1\%   &&& 2000-2100   &&  1177    &&  3.1\%   \\
        0700-0800   &&  991     &&  5.5\%   &&& 2100-2200   &&  1253    &&  1.9\%   \\
        0900-1000   &&  1075    &&  7.1\%   &&& 0000-0100   &&  674     &&  4.2\%   \\
        \bottomrule
    \end{tabular}
    \caption{Treatment effects vary by hour-of-day for both kinds of treatments. We report the three most representative cases for each treatment. Both inflow and the diverted percentage are averaged over all observations. The results reinforce the observation that adherence to Departures treatments is higher than Arrivals treatments (which we observed when not controlling for hour-of-day in~\Cref{fig:ddt_delflow}).}
    \label{tab:response-hod}
\end{table}

\subsection{Estimating Treatment Effects}

Bayesian Linear Regression yields a posterior distribution for each parameter. We focus on the autoregressive and treatment co-efficients, i.e., $\beta, \gamma_{D}, \gamma_{A}$ from~\Cref{eq:bayes-reg-spec}. Our two regressions used as outcomes the first-order rate-of-change of flow $\Delta q'$ and the standardized second-order rate-of-change of flow $\Delta q''$ respectively. 
\begin{figure*}[t]
    \captionsetup[subfigure]{justification=centering}
    \begin{subfigure}{\textwidth}
        \centering
        \includegraphics[width=0.75\textwidth]{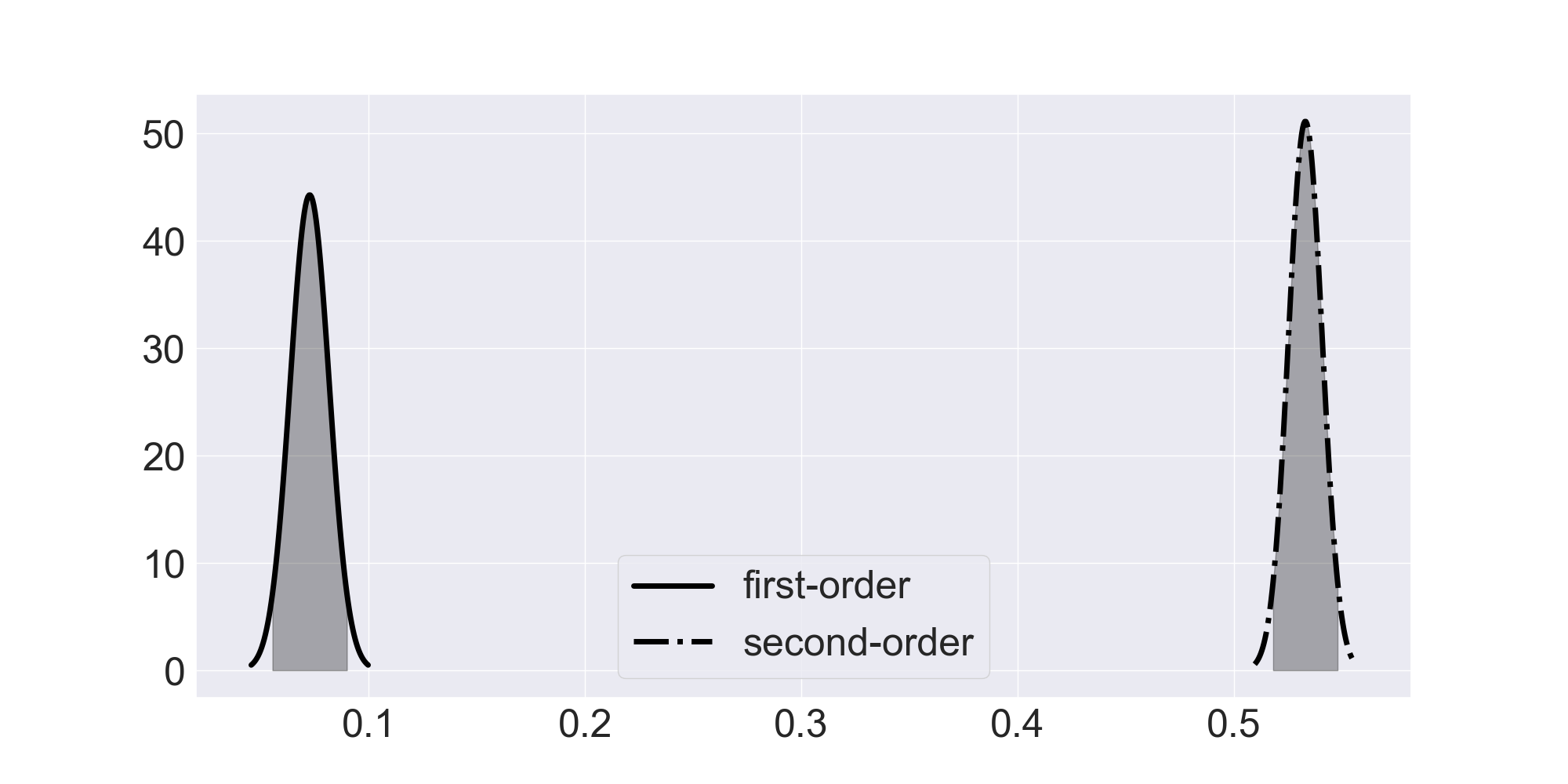}
        \caption{Autoregressive}
        \label{fig:delflow_autoreg}
    \end{subfigure}
    
    \begin{subfigure}{\textwidth}
        \centering
        \includegraphics[width=0.75\textwidth]{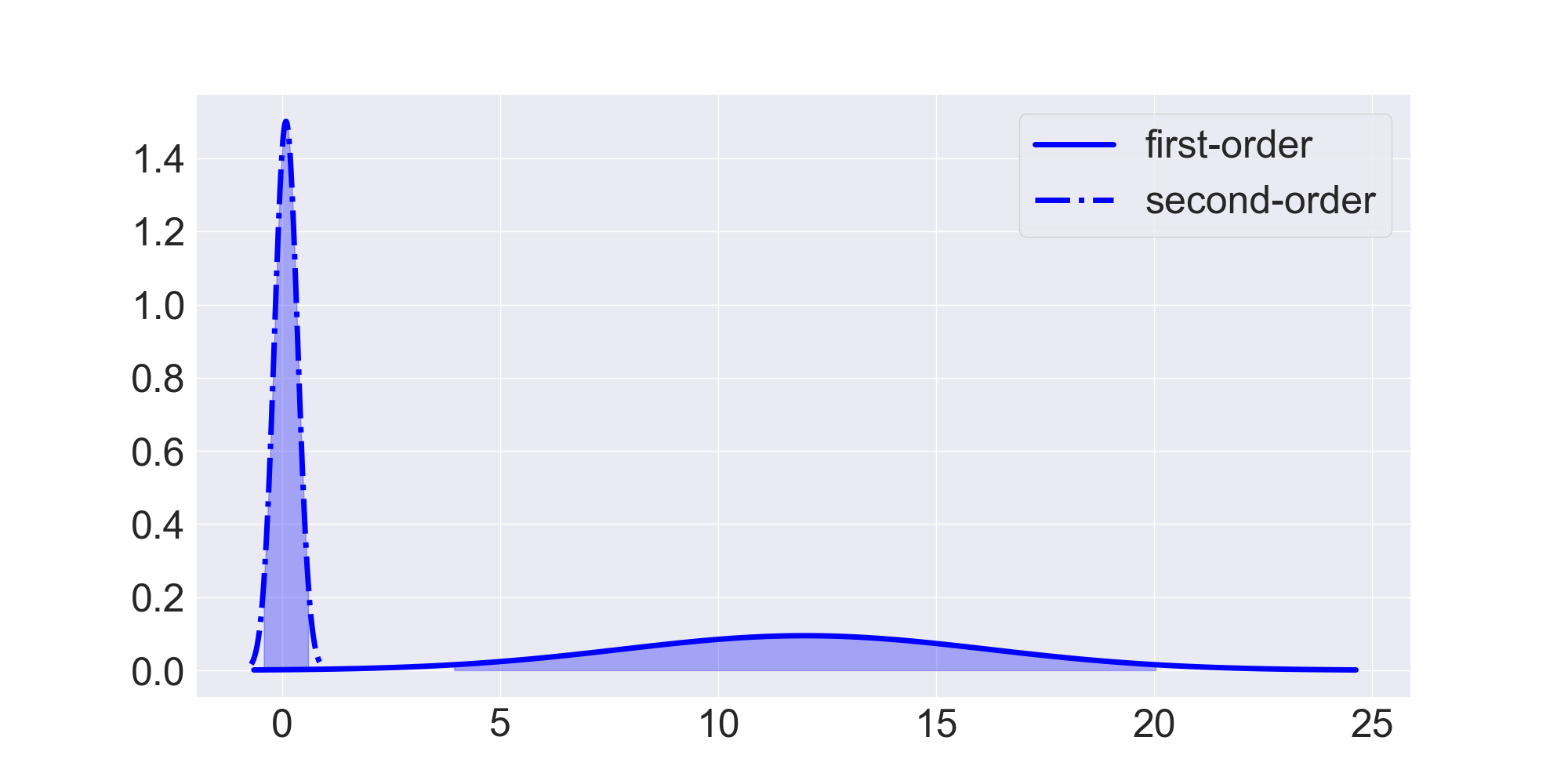}
        \caption{Departures Treatment}
        \label{fig:delflow_treatdep}
    \end{subfigure}
    
    \begin{subfigure}{\textwidth}
        \centering
        \includegraphics[width=0.75\textwidth]{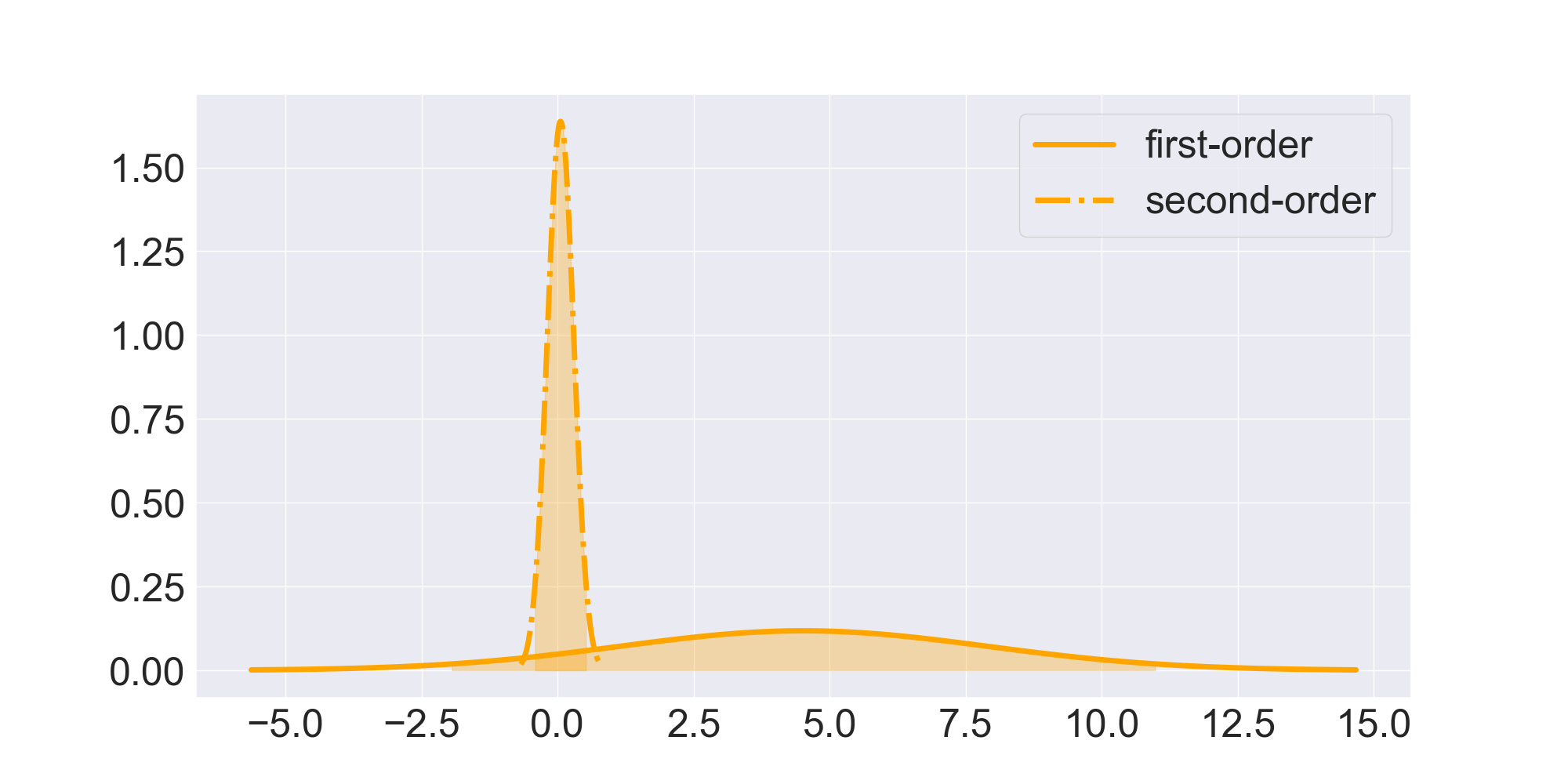}
        \caption{Arrivals Treatment}
        \label{fig:delflow_treatarr}
    \end{subfigure}
    \caption{The posterior distributions for two sets of Bayesian Linear Regression co-efficients. The \textbf{solid lines} represent the regression where the outcome is the \textbf{first-order} rate of change of difference in flow, $\Delta q'$. The auto-regressive co-efficient (\subref{fig:delflow_autoreg}) is close to 0, which implies that treatments explain most of the variation when active. The subplots (\subref{fig:delflow_treatdep}) and (\subref{fig:delflow_treatarr}) encode the average number of vehicles that respond to the Departures and Arrivals treatment respectively.
    The \textbf{dashed lines} represent the regression where the the outcome is the \textbf{standardized second-order} rate of change of difference in flow, $\Delta q''$. The autoregressive co-efficient has much higher relative magnitude than either treatment coefficient (both of which are close to 0), which suggests that the treatments do not explain much of the variation in the outcome when active.}
    \label{fig:ddt_delflow}
\end{figure*}
\Cref{fig:ddt_delflow} illustrates posterior distributions for each of $\beta, -\gamma_{D}, \gamma_{A}$ for the first-order (solid) and standardized second-order (dashed) outcomes; \textcolor{black}{we plot the distribution of $-\gamma_D$ to easily interpret it as the positive number of vehicles diverted}. The first-order auto-regressive co-efficient is close to 0, while the treatments are not, which suggests that the treatments explain most of the variation. The magnitude of the second-order auto-regressive co-efficient is much higher than either treatment co-efficient, which suggests that the treatment has no significant effect on higher-order rates-of-change.

\textit{The first-order co-efficients from \Cref{fig:ddt_delflow} suggest that on average, roughly 12 cars per 15 mins respond to Departures treatments and 4 cars per 15 minutes respond to Arrivals treatments}. Since traffic varies by hour-of-day, we expect the response to as well.  Thus, we reran the $\Delta q'$ regressions while controlling for hour-of-day. \Cref{tab:response-hod} reports the average hourly response to each treatment as a percentage of average hourly traffic, for three selected hours-of-day (roughly corresponding to AM/PM peaks). We observe considerable variation by hour, as well as clearer picture of how the response to one treatment is higher than the other. Note that response rates can change depending on several factors, such as the location of VMS board; future work will further investigate the factors that could affect the response rates.

\subsubsection*{Acknowledgements}

Pacific Northwest National Laboratory is operated by Battelle Memorial Institute for the U.S. Department of Energy under Contract No. DE-AC05-76RL01830. This work was supported by the U.S. Department of Energy Vehicle Technologies Office. We also thank the landside operations team at Seattle-Tacoma International Airport for sharing data and insights.

\bibliographystyle{apalike}
\bibliography{references}

\end{document}